\documentclass[12pt]{amsart} 

\usepackage{amsmath,amssymb,amsthm,amsfonts,enumitem}
\usepackage{color}
\usepackage{graphicx,subcaption }
\usepackage{hyperref}
\usepackage{epstopdf}
\usepackage{bm}

\theoremstyle{plain} 
\newtheorem{thm}{Theorem}[section]
\newtheorem{lem}[thm]{Lemma}
\newtheorem{prop}[thm]{Proposition}

\theoremstyle{definition}
\newtheorem{defn}{Definition}[section]
\newtheorem{conj}{Conjecture}[section]

\theoremstyle{remark}

\begin{document}

\title{On the square peg problem}

\author{Gregory R. Chambers}
\address{Department of Mathematics, Rice University, Houston, TX, 77005}
\email{gchambers@rice.edu}
\date{\today}
\begin{abstract}
	We show that if $\gamma$ is a Jordan curve in $\mathbb{R}^2$ which is close to a $C^2$ Jordan
	curve $\beta$ in $\mathbb{R}^2$, then $\gamma$ contains an inscribed square.  In particular, if $\kappa > 0$
	is the maximum unsigned curvature of $\beta$ and there is a map $f$ from the image of $\gamma$ to the image of
	$\beta$ with $||f(x) - x|| < \frac{1}{10 \kappa}$ and $f \circ \gamma$ having winding number
	$1$, then $\gamma$ has an inscribed square of positive sidelength.
\end{abstract}
\maketitle

\section{Introduction}

In 1911, O. Toeplitz made the following conjecture:
\begin{conj}[Square Peg Problem]
	\label{conj:square_peg}
	Suppose that $\gamma$ is a Jordan curve in $\mathbb{R}^2$.  There there are four distinct
	points on $\gamma$ that form a square with positive sidelength.
\end{conj}

Despite being over a century old, this conjecture still remains open.  However, inscribed squares
are known to exist for many Jordan curves.  In \cite{BM1}, B. Matschke gives an extensive survey of
the problem and known results.  In particular, in his thesis \cite{BM2}, he proved that all
Jordan curves with a technical condition (no special trapezoids of a certain size) contain inscribed squares.  More
recent articles include \cite{T}, in which T. Tao uses methods different to those discussed below
to show that certain Lipschitz Jordan curves have inscribed squares, and \cite{GL}, in which
J. E. Greene and A. Lobb show that smooth Jordan curves contain rectangles of all side ratios.

The first progress on Conjecture~\ref*{conj:square_peg} (which will also be of importance to this article) was
due to L. G. Schnirelman \cite{S}; in 1929 he proved the following theorem:
\begin{thm}
    \label{thm:std}
    Suppose that $\gamma$ is a $C^2$ Jordan curve in $\mathbb{R}^2$; then
    an arbitrarily small perturbation of $\gamma$ contains an odd number
    of squares.  In particular, every $C^2$ curve contains an inscribed square
    of positive sidelength.
\end{thm}
He actually proved somewhat more than this, and an extended and corrected version of his article
was published posthumously in 1944.  The idea is as follows.  We first note that there is an isotopy
from our curve $\gamma$ to an ellipse $\mathcal{E}$.  In addition, we can show that $\mathcal{E}$ contains
exactly one inscribed square of positive sidelength.  Next, for every curve $\beta$ in the isotopy,
we consider three manifolds $M_1, M_2,$ and $M_3$.  $M_1$ is the set of all ordered quadruples 
in $\mathbb{R}^2$, $M_2$ is the set of all ordered quadruples in $\mathbb{R}^2$ that form
the vertices of squares in counterclockwise order (including squares of $0$ sidelength),
and $M_3$ is the set of all ordered quadruples of points on $\beta$.

We see that $\beta$ has an inscribed square if $M_2$ and $M_3$ intersect away from the boundary of $M_2$, and $M_2$ and $M_3$ are both
submanifolds of $M_1$.  By perturbing the isotopy slightly, we may assume it is generic, and 
by the parametric version of Thom's Multijet Transversality Theorem (see \cite{GG}), we see that
these intermediate inscribed squares form paths from the initial inscribed squares to final inscribed squares,
or to the boundary of $M_2$.  We say squares here because each square will be counted multiple times (since
the quadruples are taken to be ordered).  If we take equivalence classes based on these equivalent squares,
then we have paths of inscribed squares whose endpoints lie on the boundary of $M_2$, correspond
to inscribed squares on the $\mathcal{E}$, or correspond to inscribed squares on $\gamma$.

Since $\gamma$ is $C^2$, we may assume that our isotopy is $C^2$, and so has universally bounded curvature.
This ensures that any inscribed square in any intermediate curve must have sidelength bounded universally below,
which in turn implies that a path of inscribed squares cannot reach the boundary of $M_2$.  Since
$\mathcal{E}$ has an odd number of squares (it has exactly 1), $\gamma$ must also have an odd number
of inscribed squares, and so cannot have zero squares.  Note that this is where this approach fails
to prove that there are is an inscribed rectangle of prescribed sidelength ratio; an ellipse has two such
rectangles, and so this approach only implies that the final curve has an even number of such rectangles,
which does not preclude it having none.

From this approach, however, we can prove the following theorem for special Jordan curves:
\begin{thm}
    \label{thm:annulus}
    Suppose that $\gamma$ is a Jordan curve which lies inside an annulus with inner radius $r$
    and outer radius $R$ with
    $$ r \leq R \leq (1 + \sqrt{2}) r.$$
    If $\gamma$ is homotopic to the outer boundary component, then it contains an inscribed square.
\end{thm}
The idea for the proof of Theorem~\ref*{thm:annulus} is to approximate $\gamma$ by a sequence
$\{ \alpha_i \}$ of $C^2$ curves which converge to $\gamma$ in $C^0$.  We then observe that since $\gamma$
is homotopic to the outer boundary component, there is a $C^2$ homotopy from the outer boundary to $\alpha_i$.
If we perturb this homotopy, we can assume that the initial curve has one square
and it has large sidelength.  Furthermore, if we follow the path of this square through the $C^2$ homotopy, the
sidelength of the square changes continuously.  In order for $\gamma$ to have no inscribed squares, for every
$\epsilon > 0$, we can find some $i$ so that $\alpha_i$ contains only inscribed squares of sidelength less
than $\epsilon$.  Choosing $\epsilon$ sufficiently small, such a scenario would force some intermediate curve to
contain a square of intermediate sidelength.  We can show that such a square would have to have
a corner inside the inner boundary component, which is impossible.

The main result of this article is a generalization of Theorem~\ref*{thm:annulus}:
\begin{thm}
	\label{thm:main}
	Suppose that $\gamma$ is a $C^2$ Jordan curve in $\mathbb{R}^2$ with unsigned curvature bounded
	above by $\kappa > 0$.  Suppose that $\beta$ is a Jordan curve
	in $\mathbb{R}^2$ so that there is a continuous function $f$ from the image of $\beta$ to the image of $\gamma$
	with the properties that
		\begin{enumerate}
			\item	$$ || f(x) - x|| < \frac{1}{10 \kappa}.$$
			\item	The map $f \circ \beta$ from $S^1$ to the image of $\gamma$ has winding number $1$ (with respect 
			to the identity function on $\gamma$).
		\end{enumerate}
	Then $\beta$ contains an inscribed square.
\end{thm}

Our technique to proving Theorem~\ref*{thm:main} is similar to the proof of Theorem~\ref*{thm:annulus}.  In particular, we approximate
$\gamma$ in $C^0$ by a sequence of $C^2$ curves $\{ \alpha_i \}$.  We then argue that there is a generic
$C^2$ homotopy from $\beta$ to each $\alpha_i$, that $\beta$ starts with an odd number of large squares,
and in the process of the homotopy only an even number of large squares are introduced or eliminated.  Thus, 
if all of the final squares in $\alpha_i$ are small, then at some point in the homotopy there must be a
square of intermediate size, which we show is impossible.

\begin{figure}
	\caption{A $C^2$ curve with small curvature and an inscribed square with short sides}
	\centerline{\includegraphics[width=0.5\textwidth]{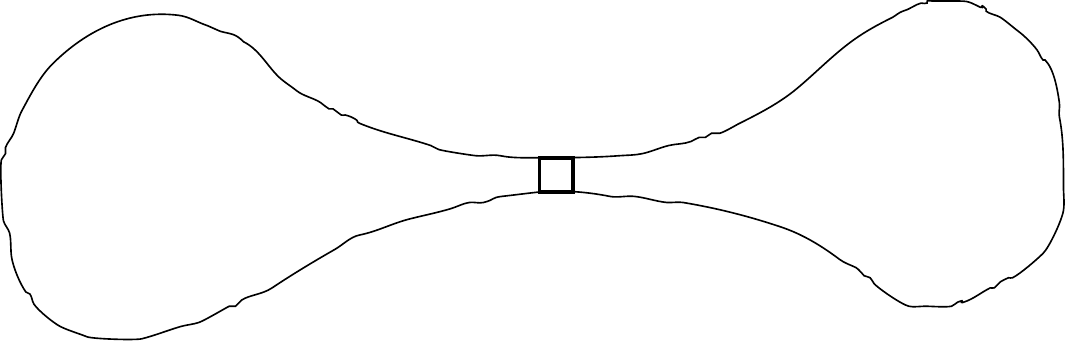}}
	\label{fig:peanut}
\end{figure}

The problem with using sidelength for ``small", ``medium", and ``large" in this context is that
we can find a $C^2$ curve $\beta$ which has arbitrarily small curvature, and which contains an arbitrarily
small square in terms of sidelength; this is shown in Figure~\ref*{fig:peanut}.  Instead, we will use a different notion 
of the size of a square; this is defined in Section~\ref*{sec:prelim}.

The remainder of the article is structured as follows.  In Section~\ref*{sec:intermediate}, we show that if $\beta$ and $\alpha$ are two $C^2$ curves which are close, then $\alpha$ cannot contain an
inscribed square of intermediate size.  We also show that inscribed squares of small (or zero) side
length have small size in terms of our new definition, and that initial squares in $\beta$ must be large.
Finally, in Section~\ref*{sec:proof}, we prove Theorem~\ref{thm:main}.

\noindent {\bf Acknowledgments} The author was supported in part by NSF
grant DMS-1906543.

\section{Preliminaries}
\label{sec:prelim}

We begin with defining a more broad notion of inscribed squares:
\begin{defn}
    \label{defn:inscribed_square}
    Suppose that $\gamma$ is a $C^2$ curve (which may have self-intersections).  A generalized inscribed square is 
    a collection of four distinct points $\nu_1$, $\nu_2$, $\nu_3$, and $\nu_4$
    on $S^1$ so that $\gamma(\nu_1)$, $\gamma(\nu_2)$, $\gamma(\nu_3)$, 
    and $\gamma(\nu_4)$ form the vertices of a square.  Note
    that this square may have zero sidelength.  For the remainder of the article,
    we will use ``inscribed square" to refer to a generalized inscribed square.
\end{defn}

We can extend Theorem~\ref*{thm:std} to include homotopies of
curves which may self-intersect.
\begin{thm}
    \label{thm:std_extended}
    Suppose that $\beta$ is a $C^2$ Jordan curve, $\alpha$ is
    a $C^2$ curve (which potentially contains self-intersections),
    and $H$ is a generic $C^2$ homotopy from $\beta$ to $\alpha$.  Then
    inscribed squares of positive sidelength are created
    and destroyed in one of the following ways:
    \begin{enumerate}
        \item Two squares merge to a square of positive side
            and are destroyed.
        \item A square of positive sidelength is created, and then
            immediately splits into two squares (of positive sidelength).
        \item A square of zero sidelength is created, and then
            grows.
        \item A square of positive sidelength shrinks to a square of zero sidelength
            and then is destroyed.
    \end{enumerate}
\end{thm}
\begin{proof}
    This uses the same cobordism techniques as are employed in
    the proof of Theorem~\ref*{thm:std} as described in the introduction.  In particular,
    the only difference is that paths of squares can hit the boundary of $M_2$; this corresponds
    exactly to squares of sidelength zero being created or destroyed.
\end{proof}

\begin{figure}
    \caption{The size of an inscribed square, with $a_i = f(a_i)$}
    \centerline{\includegraphics[width=0.5\textwidth]{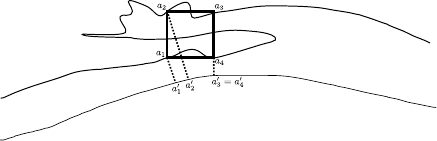}}
    \label{fig:size}
\end{figure}

We now move on to our definition of the ``size" of an inscribed square:
\begin{defn}
    \label{defn:size}
    Suppose that $\beta$ is a $C^1$ Jordan curve in $\mathbb{R}^2$ and $\gamma$ is a $C^1$ Jordan curve in $\mathbb{R}^2$.  Suppose in addition
    that $f$ is a map from the image of $\beta$ to the image of $\gamma$ so that the curve $f \circ \beta$ is a piecewise
    $C^1$ curve.  If $a_1,$ $a_2,$ $a_3,$ and $a_4$ is a sequence of vertices of a quadrilateral on $\beta$ which are ordered
    in the same order as the orientation as $\beta$, then the \emph{size} of that quadrilateral is computed as follows.
    
    We begin by dividing $\beta$ up into four segments, $\omega_1, \dots, \omega_4$ where $\omega_i$ goes from $a_i$ to $a_{i+1}$
    (with $a_5 = a_1$), and define $\tilde{\omega}_i = f \circ \omega_i$.  Every $\tilde{\omega}_i$ is a piecewise $C^1$ curve.  We
    define its oriented length, denoted by $L_o(\tilde{\omega}_i)$ as follows.  For each segment of $\tilde{\omega}_i$ which has the same
    orientation as $\gamma$, we compute its length, and for each segment of $\tilde{\omega}_i$ which has the opposite
    orientation as $\gamma$, we compute the negative of its length.  To compute $L_o(\tilde{\omega}_i)$, we add all of these quantities
    together.
    
    We now define the size of the quadrilateral with vertices $a_1, \dots, a_4$ to be
        $$ \textrm{size}(a_1,\dots,a_4) = \inf_i \sum_{j \neq i} |L_o(\tilde{\omega}_j)|.$$
    This is shown in Figure~\ref*{fig:size}.
\end{defn}

\section{No intermediate squares}

\label{sec:intermediate}

The purpose of this section is threefold.  First, we should that a $C^2$ Jordan curve $\beta$ only contains large
squares.  Second, we prove that, if $\alpha$ and $\beta$
are two $C^2$ which are close (in terms of the maximum unsigned curvature of $\beta$),
then $\alpha$ cannot have an inscribed square of intermediate size.  Third,
we show that every inscribed square in $\alpha$ which has short sidelength must have small
size.  We begin with the first result; for it, we will need several lemmas.

\begin{lem}
    \label{lem:initial_square_size}
    Suppose that $\beta$ is a $C^2$ curve with maximum unsigned curvature
    $\kappa > 0$.  Then every square in $\beta$ must have size (with respect to the identity function on $\beta$) at least
        $$ \frac{\pi}{\kappa}.$$
\end{lem}
\begin{proof}
    Let $S$ be an inscribed square of $\beta$ with consecutive vertices
    $a_1, a_2, a_3,$ and $a_4$.  Let $\ell$ be the size of $S$,
    and without loss of generality we may assume that the arc from $a_1$ to $a_2$
    to $a_3$ to $a_4$ has length $\ell$.  We observe that the tangent vectors from $a_1$
    to $a_4$ turns through at least $\pi$ radians.
    
    Since the curvature of $\beta$ is at most $\kappa$, the total angle 
    that is swept out by the tangent vectors is at most $$ \ell \kappa, $$
    and so $$\ell \geq \frac{\pi}{\kappa}.$$
\end{proof}

\begin{lem}
    \label{lem:shortest_distance}
    Suppose that $\gamma$ is a $C^2$ arc of maximum unsigned curvature
    $\kappa > 0$, and $\gamma$ has length $\ell$ at most $$\frac{\pi}{4 \kappa}.$$
    Then
        $$ ||a - b|| \geq \frac{\ell}{\sqrt{2}}, $$
    where $a$ and $b$ are the endpoints of the arc.
\end{lem}
\begin{proof}
    Suppose that $u$ and $v$ are the unit tangent vectors at $a$ and $b$; the angle between them
    is at most $\pi/4$ (due to the constraint on $\ell$).  As a result, if $v$ is the vector from $a$
    to $b$, then the projection of every tangent vector onto $v$ is of size at least $\frac{1}{\sqrt{2}}$.
    This is because there is a tangent vector on $\gamma$ which is equal to $\frac{v}{||v||}$
    (with $||v|| = ||a - b||$).
    Thus, $$ ||a - b|| \geq \int_0^\ell || \pi_v \gamma'(s)||  ds \geq \frac{\ell}{\sqrt{2}}, $$
    where $\pi_v \gamma'(s)$ is the projection of $\gamma'(s)$ onto $v$ (here, $\gamma$ has
    a unit-speed parametrization).
\end{proof}

\begin{prop}
    \label{prop:no_intermediate}
    Suppose that $\alpha$ and $\beta$ are two $C^2$ curves,
    and $\kappa > 0$ is the maximum unsigned curvature of $\beta$.
    Let $\delta = \frac{1}{10 \kappa} > 0$, and suppose that
    there is a function $f$ from the image of $\alpha$ to the image
    of $\beta$ so that 
        $$ |f(\alpha(s)) - \beta(s)| \leq \delta.$$
    Note that $\alpha$ and $\beta$ need not be unit-speed parametrizations.
    
    Then there is no inscribed square $S$ in $\alpha$ that
    has size (with respect to $f$) equal to $$\frac{\pi}{4 \kappa}.$$
\end{prop}
\begin{proof}
    Suppose that $a_1, a_2, a_3,$ and $a_4$ are the corners
    (in order) of $S$ on $\alpha$.  Let $a'_i = f(a_i)$ for all
    $i$, and assume (without loss of generality) that
    the size of $S$ is $\ell$, and is the oriented length of the
    segment from $a'_1$ to $a'_4$ (in order).
    
    Let $\theta$ and $\phi$ in $[0, \pi]$ be the angles
    at $a'_2$ and $a'_3$, respectively (with respect to
    the line segments from $a'_1$ to $a'_2$, from $a'_2$ to
    $a'_3$, and from $a'_3$ to $a'_4$.  Then the tangent vector of
    $\beta$ from $a'_1$ to $a'_4$ must pass through an angle of
    at least $$2 \pi - \theta - \phi.$$
    
    The total angle that the tangent vector goes through is at most
    $\ell \kappa$, and so
        $$ \ell \kappa \geq 2 \pi - \theta - \phi. $$
        
    Let $\mathcal{L}$ be the Euclidean sidelength of $S$.  The angle of the square $S$
    at each corner is $\pi / 2$, and when we move to the curve $\beta$, each corner can move
    a distance of at most $\delta$.  Furthermore, $\delta < \frac{1}{10 \kappa}$ and $\mathcal{L} \geq \ell / \sqrt{2} = \frac{\pi}{4 \sqrt{2} \kappa}$,
    and so $\frac{2 \delta}{\mathcal{L}} \leq 1$.  Thus, the amount that each angle can decrease is
    at most $2 \arcsin{\frac{2 \delta}{\mathcal{L}}}$, and so
        $$ 2 \pi - \theta - \phi \geq \pi - 4 \arcsin{\frac{2\delta}{\mathcal{L}}}. $$
    Suppose now that the size of $S$ is $\frac{\pi}{4 \kappa} = \ell$. Since $$\ell \leq \frac{\pi}{4 \kappa},$$
    Lemma~\ref*{lem:shortest_distance} implies that $$ \mathcal{L} \geq \frac{\ell}{\sqrt{2}}. $$
    
    The other estimate which will be relevant is that, for $x \in [0,1]$, 
        $$ \arcsin(x) \leq \frac{\pi x}{2}.$$
    
    Combining these estimates,
        $$ 1 - \frac{4 \sqrt{2} \delta}{\ell} \leq \frac{\ell \kappa}{\pi}.$$
    First, we observe that
        $$ \frac{ \ell \kappa}{\pi} = \frac{1}{4}. $$
    Next, since $$\delta < \frac{1}{10 \kappa},$$
        $$ 1 - \frac{4 \sqrt{2} \delta}{\ell} = 1 - \frac{4 \sqrt{2} \delta (4 \kappa)}{\pi} > 1 - \frac{16 \sqrt{2} }{10 \pi} > \frac{1}{4}, $$
    which cannot happen.
\end{proof}

\begin{lem}
    \label{lem:small_square}
    Suppose that $\alpha$ and $\beta$ are $C^2$ Jordan curves.  If
    $f$ is a function from the image of $\alpha$ to the image of $\beta$ so that
        $$ |f(\alpha(s)) - \beta(s) | \leq \delta < \frac{1}{10 \kappa} $$
    and $S$ is an inscribed square of $\alpha$ with size $\rho \geq 0$ with respect
    to the identity function on $\alpha$, then the size of $S$ with respect to
    $f$ is less than
        $$ \frac{\sqrt{2}}{5 \kappa} + \sqrt{2} \rho $$
\end{lem}
\begin{proof}
    As in the proof of Proposition~\ref*{prop:no_intermediate}, 
        $$ \mathcal{L} \geq \frac{\ell}{\sqrt{2}}. $$
    where $\mathcal{L}$ is the distance between the endpoints of an arc of length $\ell$ of $\alpha$.
    Since $S$ has sidelength $\rho$ and $\delta < \frac{1}{10 \kappa}$, $\mathcal{L} < \rho + \frac{1}{5 \kappa }$,
    and so
        $$ \frac{1}{5 \kappa} + \rho > \frac{\ell}{\sqrt{2}}, $$
    and so
        $$ \frac{\sqrt{2}}{5 \kappa} + \sqrt{2} \rho > \ell. $$
\end{proof}

\section{Proof of Theorem~\ref*{thm:main}}
\label{sec:proof}

In this section, we prove Theorem~\ref*{thm:main}.  We first need two lemmas.

\begin{lem}
    \label{lem:two_step_approx}
    Suppose that $\alpha$ and $\beta$ are two $C^2$ Jordan curves so that
    there is a function $f$ from the image of $\alpha$ to the image
    of $\beta$ so that
        $$ | f(\alpha(s)) - \beta(s) | < \eta $$
    for every $s \in S^1$, and so that $f(\alpha)$ has winding number
    $1$.  Then there is an is a $C^2$ isotopy $H : S^1 \times [0,1] \rightarrow \mathbb{R}^2$
    from $\beta$ to $\alpha$ along with a $C^2$ function $ P : S^1 \times [0,1] \rightarrow S^1$
    so that:
    \begin{enumerate}
        \item $|\beta(P(s,t)) - H(s,t)| < \eta$
        \item $P(s,0) = s$ for all $s \in S^1$.
        \item $f(\alpha(s)) = \beta(P(s,1))$ for all $s \in S^1$.
    \end{enumerate}
\end{lem}
\begin{proof}
    We first choose $\delta > 0$ so that $| f(\alpha(s)) - \beta(s)| < \eta - \delta$ for all
    $s \in S^1$, and so small that there is an embedded $\delta$-tubular neighborhood of
    $\beta$.  Since the map $f$ has winding number $1$, we can find a $C^2$ homotopy from $\beta$ to a $C^2$ Jordan curve $\gamma$ which lies
    in in the $\delta$-tubular neighborhood of $\gamma$, and so that
    the following holds.  If $\pi$ is the projection of the image of $\gamma$ onto
    the image of $\beta$ by projecting it in the normal direction,
    then
    $$ \pi(\gamma(s)) = f(\alpha(s)) $$
    for every $s \in S^1$.
    
    We form this isotopy by extending small loops from $\beta$ to
    form $\gamma$ (see Figure \ref*{fig:loops}).  During this portion of the homotopy,
    $P$ is defined as the projection of the images of the intermediate curves onto the image of $\beta$.
    
    \begin{figure}
	    \caption{Extending loops to produce our homotopy}
	    \centerline{\includegraphics[width=0.5\textwidth]{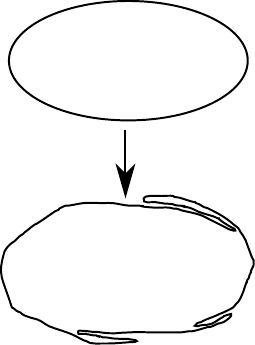}}
	    \label{fig:loops}
    \end{figure}

    To complete the homotopy, we move each element of $\gamma(s)$ to
    $\alpha(s)$ linearly (we also smooth this out so that the resulting
    homotopy is $C^2$).  For this portion of the homotopy,
    we define $P$ to take each point on each intermediate curve
    back to the same point on $\beta$ that $\gamma(s)$ was sent to.
    Since $\delta$ was chosen to be sufficiently small, the distance
    is still less than $\eta$, as desired.
\end{proof}

We can now prove Theorem~\ref*{thm:main}:
\begin{proof}[Proof of Theorem~\ref*{thm:main}]
    We begin with $\beta$.  We approximate $\gamma$ (in $C^0$)
    by a sequence $\{ \alpha_i \}$ of $C^2$ Jordan curves.
    Suppose that $\gamma$ contains no inscribed squares of
    positive sidelength.  Then, after potentially moving
    to a subsequence, we may assume the following:
    \begin{enumerate}
        \item There is a function $f_i$ from the image of $\alpha_i$ to the image of $\beta$
            so that $|f_i(\alpha_i(s)) - \beta(s)| < \eta$.  In addition,
            $f_i$ has winding number $1$.
        \item $\alpha_i$ contains an odd number of squares of positive sidelength.
        \item Each square $S$ in $\alpha_i$ has size (with respect to the identity function) at most
            $\frac{1}{i}$.
    \end{enumerate}
    Theorem~\ref*{thm:std} implies that each $\alpha_i$ has an odd number of inscribed squares of positive
    sidelength.  The function $f_i$ exists since the function $f$ exists, and $\alpha_i$ goes
    to $\gamma$ in $C^0$.  If all squares in all $\alpha_i$ for $i > N$ have size with respect to 
    $f_i$ greater than $1/N$, then we can find a sequence of squares which converge to a square $S^*$ in $\gamma$,
    and whose sizes (with respect to the identity functions) is bounded below by $1/N$.  Since $\gamma$
    contains no squares, $S^*$ must have sidelength $0$, that is, its corners must be the images of the
    same point in $S^1$.  Taking all of the $\{ \alpha_i \}$ to have unit-speed parameterizations,
    we have points $a^i_1, a^i_2, a^i_3,$ and $a^i_4$ on $S^1$ such that their images with respect
    to $\alpha_i$ form a square with size at least $1/N$ (with respect to the identity function.
    Furthermore, $a^i_1 \rightarrow a_1^*$, $a^i_2 \rightarrow a_2^*$, $a^i_3 \rightarrow a_3^*$, and
    $a^i_4 \rightarrow a_4^*$ where $a_1^*, a_2^*, a_3^*,$ and $a_4^*$ are the points which form $S^*$,
    that is, $a_1^* = a_2^* = a_3^* = a_4^*$.  Since we the parameterizations are unit speed and the size of
    each square with respect to the identity function is at least $1/N$, this convergence cannot occur.
    This is due to the fact that since $a_1^i,$ $a_2^i$, $a_3^i$, and $a_4^i$ converge to the same point,
    the length of the arc from $a_1^i$ through to $a_4^i$ must go to $0$, which does not occur.
    
    Thus, the sizes of the squares in $\alpha_i$ (with respect to the identity functions) go
    to $0$ as $i \rightarrow \infty$.  By Lemma~\ref*{lem:small_square}, the sizes of all squares in $\alpha_i$ with respect
    to $f_i$ are bounded above by $$\frac{\sqrt{2}}{5 \kappa} + \frac{1}{100 \kappa} < \frac{\pi}{4 \kappa}$$
    for $i$ sufficiently large.
    
    Let us consider now
    a homotopy from Lemma~\ref*{lem:two_step_approx} from $\beta$ to $\alpha_i$.  
    By Theorem~\ref*{thm:std_extended}, the size of squares with respect to $P$ are continuous,
    except where they are created and destroyed.
    The original curve $\beta$ contains an odd number of squares, and each has size (with respect to the identity function) at least $\frac{\pi}{\kappa}$ (by Lemma~\ref{lem:initial_square_size}).  Since all final squares have size less than $\frac{\pi}{4 \kappa}$,
    we observe that either a square along the way must have size $\frac{\pi}{4 \kappa}$, or an odd number of squares of size greater than
    $\frac{\pi}{ 4 \kappa}$ must be destroyed along the way.    By Proposition~\ref*{prop:no_intermediate}, there are no squares in any intermediate
    curve of size (with respect to $P$) equal to $\frac{\pi}{4 \kappa}$, and so the latter must be true.  Looking at Theorem~\ref*{thm:std_extended},
    since an odd number of such squares must be destroyed, at least one must be destroyed by becoming a square of $0$ sidelength and then being destroyed.  However, at the point that it becomes a square of sidelength $0$ it has size at most $\frac{\sqrt{2}}{5 \kappa} < \frac{\pi}{4 \kappa}$, and so it would have already have had to become a square of size $\frac{\pi}{4 \kappa}$ along the way, which is impossible.
    This completes the proof.
    
\end{proof}

\bibliographystyle{amsplain}
\bibliography{bibliography}

\end{document}